\documentclass[12pt]{amsart}
\usepackage{pstricks,pst-node,pst-fill,amssymb}
\usepackage[all]{xy}

\newlength{\cellsz}
\newcounter{cellsize}
\newcommand{\setcellsize}[1]{%
  \setcounter{cellsize}{#1}%
  \setlength{\cellsz}{\value{cellsize}\unitlength}}%
\setcellsize{13}%
\newcommand\cellify[1]{\def\thearg{#1}\def\nothing{}%
\ifx\thearg\nothing \vrule width0pt height\cellsz depth0pt\else
\hbox to 0pt{{\begin{picture}(\value{cellsize},\value{cellsize})
  \put(0,0){\line(1,0){\value{cellsize}}}
  \put(0,0){\line(0,1){\value{cellsize}}}
  \put(\value{cellsize},0){\line(0,1){\value{cellsize}}}
  \put(0,\value{cellsize}){\line(1,0){\value{cellsize}}} \end{picture} \hss}}\fi%
\vbox to \cellsz{ \vss \hbox to \cellsz{\hss$#1$\hss} \vss}}
\newcommand\tableau[1]{\vcenter{\vbox{\let\\\cr
\baselineskip -16000pt \lineskiplimit 16000pt \lineskip 0pt
\ialign{&\cellify{##}\cr#1\crcr}}}}
\newcommand\tabl[1]{\vtop{\let\\\cr
\baselineskip -16000pt \lineskiplimit 16000pt \lineskip 0pt
\ialign{&\cellify{##}\cr#1\crcr}}}

\title{Signed differential posets and sign-imbalance}
\author{Thomas Lam}
\address{Department of Mathematics, Harvard University, Cambridge MA, 02138, USA}
\email{tfylam@math.harvard.edu}
\urladdr{http://www.math.harvard.edu/~tfylam}
\date{November 2006}

\def\Z{{\mathbb Z}}
\def\N{{\mathbb N}}
\def\P{\mathbf{P}}

\def\L{{\mathcal L}}
\def\Fib{{\rm Fib}}
\newcommand{\ip}[1]{\langle #1 \rangle}

\newtheorem{thm}{Theorem}
\numberwithin{thm}{section}
\newtheorem{lem}[thm]{Lemma}
\newtheorem{prop}[thm]{Proposition}

\theoremstyle{definition}
\newtheorem{definition}[thm]{Definition}

\newtheorem{remark}[thm]{Remark}

\begin{document}

\begin{abstract}
We study {\it signed differential posets}, a signed version of
differential posets.  These posets satisfy enumerative
identities which are signed analogues of those satisfied by
differential posets.  Our main motivations are the sign-imbalance
identities for partition shapes originally conjectured by Stanley,
now proven in~\cite{Lam,Rei,Sjo1}.  We show that these identities
result from a signed differential poset structure on Young's
lattice, and explain similar identities for Fibonacci shapes.
\end{abstract}
\maketitle
\section{Introduction}
For a poset $P$ and a field $K$, we let $KP$ denote the vector space of finite linear combinations
of elements of $P$.
Differential posets are posets $P$ naturally equipped with linear operators
$U,D: KP \to KP$ satisfying $DU - UD = I$.  Differential posets were introduced by Stanley in~\cite{Sta88} and independently discovered by Fomin~\cite{Fom88} who called them $Y$-graphs.  Differential posets satisfy many enumerative
properties generalizing properties of Young's lattice $Y$.  These include
\begin{eqnarray}
\label{eq:nfact} \sum_{x \in P_n} f(x)^2 & =& n! \\
\label{eq:nsum} \sum_{x \in P_n} f(x) &=& \#\{w \in S_n \mid w^2 =
1\}
\end{eqnarray}
where $f(x)$ denotes the number of maximal chains from the minimum
element of $P$ to $x \in P$.

In this paper we study {\it signed differential posets}.  A signing
$(P,s,v)$ of a poset $P$ is an assignment of a sign $v(x) \in \{\pm
1 \}$ to each element $x \in P$, and a sign $s(x \lessdot y)$ to
each cover relation $x \lessdot y$ of $P$.  Given a signed poset
$P$, one defines linear operators $U,D: KP \to KP$.  Signed
differential posets are those signed posets which give rise to the
relation $DU + UD = I$. Signed differential posets come in a number
of variations, and as the most interesting example, $\beta$-signed
differential poset satisfy the enumerative identities
\begin{eqnarray}
\label{eq:signfact}\sum_{x \in P_n} v(x)e(x)^2 & =& 0 \\
\label{eq:signsum}\sum_{x \in P_n} e(x) &=& 2^{\lfloor n/2 \rfloor}
\end{eqnarray}
where $e(x)$ is a signed sum of chains in $P$, defined using the signs $s(x \lessdot y)$ of the cover relations.

\smallskip
Our investigations were motivated by identities involving the {\it
sign-imbalance} of partition shapes, a topic studied
in~\cite{Lam,Rei,Sjo1,Sta,Whi}. For a poset $P$ and a labeling
$\omega:P \to \{1,2,\ldots,|P|\}$ one can define its sign-imbalance
$I_{P,\omega} \in \Z$, as a sum of signs over all linear extensions
of $P$.  In the case that $P$ is a Young diagram $\lambda$, the
number $I_\lambda$ is the sum of the signs $s(T) \in \{\pm 1 \}$ of
the reading words of the standard tableaux $T$ of shape $\lambda$.
Stanley~\cite{Sta} conjectured that the sign-imbalances $I_\lambda$
satisfy identities similar to (\ref{eq:signfact}) and
(\ref{eq:signsum}), with $I_\lambda$ replacing $e(\lambda)$ (and
Young's lattice $Y$ taking the place of $P$).  These (and more
general) identities were proved in~\cite{Lam, Rei, Sjo1,Sjo2}. We
show that $Y$ can be given the structure of a signed differential
poset and that the sign-imbalance identities are a consequence of
(\ref{eq:signfact}) and (\ref{eq:signsum}) which hold for all signed
differential posets.  This exhibits (\ref{eq:signfact}) and
(\ref{eq:signsum}) vividly as signed analogues of (\ref{eq:nfact})
and (\ref{eq:nsum}).

\smallskip

In Section~\ref{sec:diff}, we define {\it $\alpha$-signed differential posets} and {\it $\beta$-signed differential posets} both of which are contained in the larger class of {\it weakly signed differential posets}.  In Section~\ref{sec:id}, we give identities involving signed chains and signed walks in these classes of posets.  Our first main aim here is to generalize as many of the results in Stanley's work~\cite{Sta88} as possible.  Because of the potential for cancelation in enumerative problems for signed
differential posets, our identities are often simpler, and enumerative constants more likely to be zero, than for differential posets.

In Section~\ref{sec:examples}, we discuss our two main examples of signed differential posets: the signed Young lattice and the signed Fibonacci differential poset.  The signed Fibonacci differential
poset admits a simple construction using a signed analogue of the reflection extension which constructs
the Fibonacci differential poset.  As a consequence one can obtain in this way a large family of
signed differential posets.  It is curious that the
underlying posets for our main examples are themselves differential posets.  We have no simple
explanation of this phenomenon and have not yet found a natural signed differential poset which
is not also a differential poset.  In Section~\ref{sec:sign}, we relate signed differential posets to sign imbalance.  Besides the case of partition shapes, we show that the elegant sign-imbalance identities
are also satisfied for Fibonacci shapes.  For the case of Fibonacci shapes, one can also define sign-imbalance as the sum of signs of reading words of Fibonacci tableaux.

\smallskip
There are a number of generalizations of this work which we have not
included here to keep the connection with sign-imbalance
transparent, but we intend to investigate these in a sequel.  These
include the study of the relation $DU + UD = rI$ for $r > 1$, the
study of the $q$-analogue $DU -qUD = rI$ of all these commutation
relations, and the extension to the ``weighted'' situation, for
example in the sense of Fomin's dual-graded graphs~\cite{Fom1,Fom2}.
Some of our identities are formal consequences of Fomin's work (and
many are not), though strictly speaking Fomin disallows graphs with
negatively weighted edges, thus excluding signed differential
posets.

\smallskip

{\bf Acknowledgements.} I am grateful to Sergey Fomin for comments on an earlier version of this article.

\section{Signed differential posets}
\label{sec:diff}
\subsection{Differential posets} Let $P = \cup_{n \geq 0} P_n$ be a
graded poset with finitely many elements of each rank and with a
minimum $\hat 0 \in P_0$. We denote the partial order on $P$ by
``$<$'' and we let $x \lessdot y$ denote a cover in $P$. Thus, $x
\lessdot y$ if and only if $x < y$ and for any $z \in P$ satisfying
$x \leq z \leq y$ we have $z \in \{x,y\}$.

Let $K$ be a field of characteristic 0.  Let $KP$ denote
the $K$-vector space with basis $P$ and let $\hat KP$ denote the
$K$-vector space of arbitrary linear combinations of elements of
$P$.

A linear transformation $T: \hat KP \to \hat KP$ is called {\it
continuous} if it preserves infinite linear combinations.  Define
two continuous linear transformations $\bar U,\bar D: \hat KP \to
\hat KP$ by
\begin{align*}
\bar U x = \sum_{y \gtrdot x} y \quad \quad
\bar D x  = \sum_{y \lessdot x} y. 
\end{align*}

We use the notation $\bar U$ and $\bar D$ here instead of $U$ and
$D$ so as to save the latter notation for the signed case.  The
poset $P$ is called {\it $r$-differential} if $\bar D \bar U-\bar
U\bar D = rI$, where $I$ is the identity transformation and $r \in
\N$ is a positive integer. Differential posets were introduced by
Stanley in~\cite{Sta88} and independently studied from a more combinatorial perspective
by Fomin~\cite{Fom88} (for $r = 1$).  Our approach imitates the linear-algebraic approach used by Stanley though
many of our results can also be established in the framework of Fomin's growth diagrams.  For the purposes of the current article we
will always assume that $r = 1$, and omit the mention of $r$.

If $S \subset P$, we write $\mathbf{S} = \sum_{x\in S} x \in \hat
KP$.  The following theorem follows quickly from the definition.

\begin{thm}[{\cite[Theorem 2.3]{Sta88}}]
If $P$ is a differential poset, then
$$
\bar D\mathbf{P} = (\bar U+I)\mathbf{P}.
$$
\end{thm}

The following result are special cases of enumerative properties of
differential posets.

\begin{thm}[\cite{Fom88,Sta88}]
\label{thm:nfact} Suppose $P$ is a differential poset.
Then the identities (\ref{eq:nfact})
and (\ref{eq:nsum}) hold for each $n \in \N$.
\end{thm}

\subsection{Signed posets} Let $P = \cup_{n \in \Z} P_n$ be a graded
poset.  Denote by $E(P)$ the set of edges of its Hasse diagram, or
equivalently the set $\{x \lessdot y \mid x,y \in P\}$ of covers in
$P$.  A {\it signing} of $P$ is a pair $(s,v)$ where $s: E(P) \to
\{\pm 1\}$ is a labeling of the edges of the Hasse diagram of $P$
with a sign $\pm 1$, and $v: P \to \{\pm 1\}$ is a labeling of the
elements of $P$ with a sign.  We will call the triple $(P,s,v)$ a
{\it signed poset}.  If the signing is understood, then we may just
say that $P$ is a signed poset.

Let $(P,s,v)$ be a signed poset where $P$ has finitely many elements
of each rank. We define two continuous linear operators $U, D: \hat
KP \to \hat KP$ by
\begin{align*}
Ux &= \sum_{x \lessdot y} s(x \lessdot y) \,y \\
Dx &= \sum_{y \lessdot x} s(y \lessdot x)  \, v(x) \, v(y) \,y.
\end{align*}
 The
function $s': E(P) \rightarrow \{\pm 1\}$ given by $s'(x \lessdot y) =
 s(x \lessdot y) \, v(x) \, v(y)$ is called the {\it conjugate} of
$s$ by $v$.  The reader concerned with the asymmetry in the
definitions of $U$ and $D$ above should note that the roles are
swapped if $s$ is replaced by its conjugate $s'$.

The vector space $KP$ is naturally equipped with a symmetric
bilinear inner product $\ip{.,.}$ defined by $\ip{x,y} =
\delta_{xy}$ for $x,y \in P$. Define a deformed inner product
$\ip{x,y}_v$ by $\ip{x,y}_v = \delta_{xy}v(x)$.  Then the operators
$U$ and $D$ defined above are adjoint with respect to $\ip{x,y}_v$.

\subsection{Signed differential posets} Recall that if $S \subset P$,
we write $\mathbf{S} = \sum_{x\in S} x \in \hat KP$.

\begin{definition}
Let $(P,s,v)$ be a signed poset, where $P$ is a graded poset with a
minimum $\hat 0 \in P_0$ which has finitely many elements of each
rank. Then $P$ is {\it weakly signed differential} if we have
$v(\hat 0) = 1$ and
$$
UD+ DU = I.
$$
If, in addition, we have
\begin{equation}
\label{eq:alpha} (U + D) \P = \P
\end{equation}
then we say that $P$ is a {\it $\alpha$-signed differential poset}.
If instead we have
\begin{equation}
\label{eq:beta} (D - U) \P = \P
\end{equation}
then we say that $P$ is a {\it $\beta$-signed differential poset}.
\end{definition}

Thus a $\alpha$ or $\beta$-signed differential poset is
automatically weakly signed differential.  In the case of differential posets or in the case of (some of) Fomin's self-dual
graphs, the analogue of the equation $(U + D) \P = \P$ can be
deduced from the analogue of $UD + DU = I$.  However, in our case
such an equation is not automatically satisfied, and indeed we have
interesting examples satisfying the $\beta$-equation
(\ref{eq:beta}). Note that the relations for a $\beta$-signed
differential poset is a relation not even formally considered by
Fomin~\cite{Fom1}.

Define $\varepsilon: \N \to \{0, 1\}$ by $\varepsilon(i) = \frac{1 +
(-1)^{i+1}}{2}$. Alternatively, $\varepsilon(i) = 0$ if $i$ is even
and $\varepsilon(i) = 1$ if $i$ is odd.  The following lemma follows
formally from calculations in~\cite{Fom1}, though strictly speaking
Fomin only allows ``positive'' edges, thus excluding the relation
$UD+DU = I$.

\begin{lem}
\label{lem:signcom} Let $(P,s,v)$ be a weakly signed differential
poset. Then for each $k, l \in \N$ the linear operators $U, D$
satisfy the following relations:
\begin{align*}
DU^k &= \varepsilon(k)\, U^{k-1} + (-1)^k U^k D \\
D^kU &= \varepsilon(k)\, D^{k-1} + (-1)^k U D^k \\
D^l U^{k+l} &= U^k \prod_{s = k+1}^{k+l}((-1)^sUD + \varepsilon(s))
\\
D^{k+l} U^{l} &= \prod_{s = k+1}^{k+l}((-1)^sUD +
\varepsilon(s))D^k.
\end{align*}
\end{lem}
Note that the factors in the above products commute, so their order
is not important.

\begin{proof}
The first and second equations follows by induction from $DU = I -
UD$. The third equation follows from the first by induction from the
calculation
$$
D((-1)^sUD + \varepsilon(s)) = 
((-1)^{s+1}UD + \varepsilon(s+1))D.
$$
The last equation follows from a similar calculation.
\end{proof}
Note that we may swap $U$ and $D$ in any identity we deduce from the
identity $UD + DU = I$.  If $w \in \{U,D\}^*$ is a word in the
letters $U$ and $D$, we let $\bar w$ be $w$ reversed.  Since $UD +
DU = I$ is $\bar{ }\;$-invariant, we may also reverse the order of
all words in the identities of Lemma~\ref{lem:signcom}.

\begin{lem}
\label{lem:skew} Let $(P,s,v)$ be a weakly signed differential poset
and $n \in \N$. Then \begin{align*} D^nU^n = \begin{cases} U^nD^n &
\mbox{if $n$ is even,} \\ U^{n-1}D^{n-1} - U^nD^n & \mbox{if $n$ is
odd.}
\end{cases}
\end{align*}
\end{lem}
\begin{proof}
We proceed by induction.  The result follows from $UD + DU = I$ when
$n = 1$.  Now suppose $n$ is even and the result has been shown for
smaller $n$.  Then using Lemma~\ref{lem:signcom} repeatedly,
\begin{align*}
D^nU^n &= D(U^{n-2}D^{n-2} - U^{n-1}D^{n-1})U \\
&= U^{n-2}DD^{n-2}U - (U^{n-2} - U^{n-1}D)D^{n-1}U \\
&= U^nD^n.
\end{align*}
Now suppose that $n$ is odd.  We calculate using the inductive
hypothesis
$$
D^nU^n = D(U^{n-1}D^{n-1})U = U^{n-1}D^n U = U^{n-1}D^{n-1}-U^nD^n.
$$
\end{proof}

\section{Signed chains in signed differential posets}
\label{sec:id}
\subsection{Enumeration for weakly signed differential posets} We
collect here a few enumerative results which can be proved for all
weakly signed differential posets.

Let $(P,s,v)$ be a weakly signed differential poset.  For $x \in
P_r$, let $S(x) = \{C = (\hat 0 \lessdot x_1 \lessdot x_2 \lessdot
\cdots \lessdot x_r = x)\}$ denote the set of maximal chains $C$
from $\hat 0$ to $x$.  Define the signed sum of chains
$$
e(x) = \sum_{C =(x_i) \in S(x)} s(\hat 0 \lessdot x_1) s(x_1
\lessdot x_2) \cdots s(x_{r-1} \lessdot x_r).
$$

The following lemma is immediate from the definitions.
\begin{lem}
\label{lem:chain} For $x \in P_r$ we have
$$e(x) = \ip{U^r \, \hat 0, x} = v(x)\ip{U^r \,\hat 0, x}_v = v(x) \ip{D^r\, x, \hat 0}_v = v(x) \ip{D^r\, x, \hat 0}.$$
\end{lem}

The following result is the signed analogue of (\ref{eq:nfact}).
\begin{thm}
\label{thm:signnfact} Let $(P,s,v)$ be a weakly signed differential
poset. Then for $n \geq 2$,
$$
\sum_{x \in P_n} v(x)e(x)^2 = 0.
$$
\end{thm}
\begin{proof}
By Lemma~\ref{lem:chain},
$$\ip{D^nU^n \, \hat 0, \hat 0} =\sum_{x \in P_n}\ip{U^n \, \hat 0, x}\ip{D^n \,x , \hat 0}
=\sum_{x \in P_n} v(x)e(x)^2.$$ By Lemma~\ref{lem:signcom}, we have
$D^nU^n = \prod_{s = n+1}^{2n}((-1)^s UD + \varepsilon(s))$. Since
$\varepsilon(s) = 0$ when $s$ is even and $UD \, \hat 0 = 0$ we have
$\ip{D^nU^n \, \hat 0, \hat 0} = 0$ for $n \geq 2$.  Alternatively,
we could have used Lemma~\ref{lem:skew}.
\end{proof}

Theorem~\ref{thm:signnfact} can be generalized to show that other
enumerative invariants of a weakly signed differential poset $P$ are
independent of the poset $P$.  We now give two examples of this.

Suppose that $w \in \{U, D\}^*$ is a word in the letters $\{U,D\}$.
We say that $w$ has {\it rank} $\rho(w) = l$ if the difference
between the number of $U$'s and the number of $D$'s in $w$ is equal
to $l$.  We say that $w$ {\it vanishes} if there is a representation
of $w$ as a concatenation $u\,D\,v$ for $u,v \in \{U, D\}^*$ such
that $\rho(v)$ is even.

\begin{thm}
\label{thm:w}
Let $(P,s,v)$ be a weakly signed differential poset.  Suppose that
$w \in \{U,D\}^*$ and $x \in P$ satisfy $\rho(w) =\rho(x)$.  Then
$$
\ip{w \, \hat 0,x} = \begin{cases} 0 & \mbox{$w$ vanishes,} \\
                                e(x) & \mbox{otherwise.}
                                \end{cases}
$$  In other words, if $e(x) \neq 0$, we have $\frac{\ip{w \, \hat 0,x}}{e(x)} \in \{0,1\}$ not
depending on $P$.
\end{thm}
\begin{proof}
For each word $w \in \{U,D\}^*$ one can write, using the relation
$DU + UD = I$ only, $w = \sum_{i,j} c_{ij}(w) U^iD^j$ for some
coefficients $c_{ij}(w) \in \Z $.  The coefficient $c_{ij}(w)$ is
zero unless $i - j = \rho(w)$.  We explain why the $c_{ij}(w)$ are
unique later in Remark~\ref{rem:unique}.  For now we
imitate~\cite{Sta88} and give a method to calculate the $c_{ij}(w)$
unambiguously.  It is easy to see that $c_{ij}(Uw) = c_{i-1,j}(w)$.
We also have by Lemma~\ref{lem:signcom},
\begin{align*}
Dw &= \sum_{i,j} c_{ij}(w) DU^iD^j \\
&= \sum_{i,j} c_{ij}(w) (\varepsilon(i)U^{i-1} + (-1)^{i}U^iD)D^j.
\end{align*}
Thus $c_{ij}(Dw) = \varepsilon(i+1)c_{i+1,j}(w) +(-1)^ic_{i,j-1}$.
 When $j = 0$, we have
\begin{align*}
c_{i,0}(Uw) &= c_{i-1,0}(w) \\
c_{i,0}(Dw) &= \varepsilon(i+1)c_{i+1,0}(w).
\end{align*}
Thus $c_{\rho(w),0}(w) = 0$ or $1$ depending on whether $w$
vanishes. Since $D^j \hat 0 = 0$ for $j > 0$, we have
$$\ip{w \,\hat 0, x} = c_{\rho(w),0}\,\ip{U^{\rho(w)}\; \hat 0,x} =
c_{\rho(w),0}\,e(x).$$
\end{proof}

\begin{remark}
The coefficients $c_{ij}(w)$ are the {\it signed normal order coefficients}.  The normal
order coefficients have many interpretations, for example as rook numbers on a Ferrers board.  The signed normal
order coefficients can be obtained as specializations of the $q$-normal order coefficients~\cite{Var}.  The $q$-normal order
coefficients correspond to the relation $DU - qUD = I$ which we will study in a separate article.
\end{remark}

Define the rank generating function $F(P,t)= \sum_{x \in P}
t^\rho(x)$, where $\rho:P \to \N$ is the rank function of $P$. Let
$k,n \in \N$ and define integers $\kappa_{n,k} = \sum_{x \in P_n}
\ip{D^kU^k\,x,x}$. Now define $F_k(P,t) = \sum_{n \geq 0}
\kappa_{n,k}t^n$.  Thus $F(P,t) = F_0(P,t)$.

In~\cite{Sta88}, Stanley proved in the case of (non-signed)
differential posets that the ratios $F_k(P,t)/F(P,t)$ were rational
functions not depending on $P$.  Here we obtain the signed-analogue
of this result.

\begin{thm}
\label{thm:F} Let $(P,s,v)$ be a weakly signed differential poset.
Then
$$
F_k(P,t)= \begin{cases} F(P,t) & \mbox{if $k = 0$,} \\
F(P,t)/(1+t) & \mbox{if $k = 1$,} \\
0 & \mbox{if $k \geq 2$.} \\
\end{cases}
$$
\end{thm}
\begin{proof}
Note that $\kappa_{n,k} = \sum_{x \in P_n} \ip{D^kU^k\,x,x} =
\sum_{x \in P_{n+k}} \ip{U^kD^k\,x,x}$.  Using Lemma~\ref{lem:skew},
we calculate,
\begin{align*}
\kappa_{n,k} &= \sum_{x \in P_n} \ip{D^kU^k\,x,x} \\
&= \begin{cases} \sum_{x \in P_n} \ip{U^kD^k \, x,x} =
\kappa_{n-k,k}& \mbox{if $k$ is even, } \\ \sum_{x \in P_n}
\ip{(U^{k-1}D^{k-1}-U^kD^k) \, x,x} = \kappa_{n-k+1,k-1} -
\kappa_{n-k,k}& \mbox{if $k$ is odd.}
\end{cases}
\end{align*}
Thus,
\begin{align*}
F_k(P,t) &= \sum_{n \geq 0} \kappa_{n,k}t^n \\
&= \begin{cases} \sum_{n \geq 0}\kappa_{n-k,k}t^n = t^kF_k(P,t) & \mbox{if $k$ is even, } \\
\sum_{n \geq 0}(\kappa_{n-k+1,k-1} - \kappa_{n-k,k})t^n =
t^{k-1}F_{k-1}(P,t) - t^kF_k(P,t)& \mbox{if $k$ is odd.}
\end{cases}
\end{align*}
So we have $F_k(P,t) = 0$ if $k > 0$ is even and by definition
$F_0(P,t) = F(P,t)$. We also have
$$
F_k(P,t) = \frac{t^{k-1}F_{k-1}(P,t)}{1+t^k}
$$
if $k$ is odd, giving the stated result.
\end{proof}


\subsection{Enumeration for $\alpha$ and $\beta$-signed differential
posets} In addition to the enumerative properties shared by all
weakly signed differential posets, the $\alpha$ and $\beta$-signed
differential posets satisfy more enumerative identities.

Using the relations of a signed differential poset, it is easy to
see that there are polynomials $g_k^\alpha(z),g_k^\beta(z) \in
\Z[z]$ such that $D^k \, \P = g_k^\alpha(U) \,\P$ and $D^k \, \P =
g_k^\beta(U) \,\P$ in all $\alpha$- or $\beta$-signed differential
posets.  We will explain in Remark~\ref{rem:unique} why these
polynomials are unique. (For now, one may think of them as defined
modulo the ideal $I = \{f(z) \mid f(U)\,\P =0\} \subset \Z[z]$.)


\begin{lem}
\label{lem:alpha} Suppose $(P,s,v)$ is a $\alpha$-signed
differential poset. Then for $k \in \N$ and $i \in \{1,2,3,4\}$, we
have
\begin{align*}
g_{4k+1}^\alpha(z) & = z^{4k} -z^{4k+1} \\
g_{4k+2}^\alpha(z)& = -z^{4k+2} \\
g_{4k+3}^\alpha(z)& = -z^{4k+2} +z^{4k+3} \\
g_{4k+4}^\alpha(z)& = z^{4k+4}. \\
\end{align*}
\end{lem}
\begin{proof}

Using Lemma~\ref{lem:signcom} and (\ref{eq:alpha}), we have
$$
DU^k\,\P = (\varepsilon(k)U^{k-1} + (-1)^kU^k + (-1)^{k+1}U^{k+1})
\P,
$$
and the result follows from induction with a case-by-case analysis.
\end{proof}

\begin{lem}
\label{lem:beta} Suppose $(P,s,v)$ is a $\beta$-signed differential
poset.  Then for $l \in \N$, we have
$$
g^{\beta}_{2l}(z) = (2 - z^2)^l.
$$
The polynomial $g^\beta_{2l+1}(z)$ can be obtained from
$g^\beta_{2l}(z)$ by applying the linear transformation $z^{2i}
\mapsto z^{2i}+z^{2i+1}$ on $\Z[z]$.
\end{lem}
\begin{proof}
Using Lemma~\ref{lem:signcom} and (\ref{eq:alpha}), we have
\begin{equation}
\label{eq:DUbeta} DU^k\,\P = (\varepsilon(k)U^{k-1} + (-1)^kU^k +
(-1)^{k}U^{k+1}) \P,
\end{equation}
for any $k \in \N$.  The second statement of the theorem does follow
easily from the first.  Iterating (\ref{eq:DUbeta}) twice for $k$
even we have
$$
D^2U^{2l} \, \P = (2U^{2l} - U^{2l+2}) \P,
$$
proving the first statement by induction.
\end{proof}

Our first theorem here is the signed analogue of (\ref{eq:nsum}).
\begin{thm}
\label{thm:signnsum} Let $(P,s,v)$ be a $a$-signed differential
poset where $a \in \{\alpha,\beta\}$.  Then for $n \geq 2$, we have
\begin{align*}
\sum_{x \in P_n} v(x)e(x) = \begin{cases} 0 & \mbox{if $a = \alpha$,} \\
                                    2^{\lfloor n/2 \rfloor} & \mbox{if $a = \beta$.}
\end{cases}
\end{align*}
\end{thm}

\begin{proof}
Let $\P_n = \sum_{x \in P_n} x \in KP$.  Then
$$
\sum_{x \in P_n} v(x)e(x) = \ip{D^n\P_n,\hat 0}_v = \ip{D^n \P, \hat
0}.
$$
The result thus follows from Lemmas~\ref{lem:alpha} and
\ref{lem:beta}, since $\ip{U^j \P, \hat 0} = 0$ for $j > 0$.
\end{proof}

Generalizing Theorem~\ref{thm:signnsum}, we define for each $k,n \in
\N$, the sum $\tau_{k,n} = \ip{D^k \P_{n+k}, \P_n}_v$.  Also let
$G_k(P,t) = \sum_{n \geq 0} \tau_{k,n}t^n$.  If we let $G(P,t) =
\sum_{x \in P} v(x)t^{\rho(x)}$ denote the $v$-weighted rank
generating function of $P$, then we have $G(P,t) = G_0(P,t)$.  The
following result is a signed analogue of~\cite[Theorem 3.2]{Sta88}.

\begin{thm}
\label{thm:G} Let $(P,s,v)$ be a $a$-signed differential poset where
$a \in \{\alpha,\beta\}$ and $k\in \N$.  Then the ratio
$G_k(P,t)/G(P,t)$ is a rational function of $t$ only depending on
$k$ and $a$.
\end{thm}
\begin{proof}
Using Lemmas~\ref{lem:alpha} and \ref{lem:beta} we can write
$g^a_k(z) = a_kz^k + \cdots + a_0$. We have
\begin{align*}
\tau_{k,n} &= \ip{D^k \P_{n+k}, \P_n}_v = \ip{D^k \P,\P_n}_v \\
&= \ip{g^a_k(U) \P,\P_n}_v = \sum_{i=0}^k a_i \ip{U^i \P,\P_n}_v \\
&= \sum_{i=0}^k a_i \ip{D^i \P_n,\P}_v = \sum_{i=0}^k a_i \ip{D^i
\P_n,\P_{n-i}}_v.
\end{align*}
Thus
\begin{align*}
G_k(P,t) &= \sum_{n \geq 0} \tau_{k,n}t^n = \sum_{i=0}^k a_i \sum_{n \geq 0} \ip{D^i \P_n,\P_{n-i}}_v t^n \\
&= \sum_{i=0}^k a_i \sum_{n \geq 0} \tau_{i,n-i}t^n = \sum_{i=0}^k
a_i t^i G_i(P,t).
\end{align*}
We may assume by induction that $G_i(P,t)$ has the form given in the
theorem for $0 \leq i < k$, and so we can rearrange to write
$G_k(P,t)$ as a rational function times $G(P,t)$.  The constants
$\{a_i\}$ do not depend on $(P,s,v)$ so we are done.
\end{proof}

For $a \in \{\alpha,\beta\}$, we denote by $A^a_k(t) =
G_k(P,t)/G(P,t)$ the rational function defined by
Theorem~\ref{thm:G}.  We can calculate $A^\alpha_k(t)$ explicitly.

\begin{prop}
\label{prop:Galpha}
Let $k \in \N$.  Then
$$
A^\alpha_k(t) = \begin{cases} 1 & \mbox{if $k = 0$,} \\
                        \frac{1}{1+t} & \mbox{if $k = 1$,} \\
                        0 &\mbox{if $k > 1 $.}
                        \end{cases}
$$
\end{prop}
\begin{proof}
The result follows immediately from the recursion in the proof of Theorem~\ref{thm:G} and Lemma~\ref{lem:alpha}.
\end{proof}

In the $\beta$-case, the polynomials $A^\beta_{k}(t)$ for even $k$
have a simple form.  The odd case appears to be considerably more complicated.

\begin{prop}
\label{prop:Gbeta} Let $l \in \N$.  Then $$A^\beta_{2l}(t)= \left(\frac{2}{1+t^2}\right)^l.$$
\end{prop}
\begin{proof}
We proceed by induction.  By definition, $A^\beta_0(t) = 1$.  Using the recursion in the proof of Theorem~\ref{thm:G} and Lemma~\ref{lem:beta} one
needs to check that $A^\beta_{2k} = (\frac{2}{1+t^2})^k$ satisfies the equality
$$
A^\beta_{2k} = \sum_{i=0}^k 2^{k-i} {n \choose i} (-1)^{i} t^{2i}\, A^\beta_{2i}.
$$
The left hand side is equal to
$$
2^k \, \sum_{i = 0}^k {n \choose i} \left(\frac{-t^2}{1+t^2}\right)^i = 2^k \left(1 - \frac{t^2}{1+t^2}\right)^k,
$$
consistent with the claimed formula.
\end{proof}

\section{Two fundamental examples}
\label{sec:examples} Our two examples of signed differential posets
come from signings of the two fundamental examples of differential
posets. While it is possible to construct trivial examples of
(weakly) signed differential posets, Theorem~\ref{thm:signnsum}
shows that a $\beta$-signed differential poset must be infinite and
non-trivial.

\subsection{Young's Lattice} Let $Y$ denote Young's lattice.  Thus $Y$
is the poset of partitions $\lambda = (\lambda_1 \geq \lambda_2 \geq
\dots \geq \lambda_l > 0)$ ordered by inclusion of Young diagrams
(see for example~\cite{EC2}).  We will often identify a partition
with its Young diagram without comment, and will always think of
Young diagrams in the English notation (top-left justified). The
rank function $\rho: Y \to \N$ is given by $\rho(\lambda) =
|\lambda| = \lambda_1 + \cdots + \lambda_l$. Thus the Young diagram
of $\lambda$ has $\rho(\lambda)$ boxes.  A partition $\mu$ covers
$\lambda$ in $Y$ if $\mu$ and $\lambda$ differ by a box.  If
$\lambda$ is a partition then $\lambda'$ denotes the conjugate
partition, obtained by reflecting the Young diagram along the main
diagonal.  Recall that an {\it outer corner} of $\lambda$ is a box
that can be added to $\lambda$ to obtain a partition, while an {\it
inner corner} is a box that can be similarly removed.

Define \begin{align*} a(\lambda) &= (-1)^{\lambda_2 + \lambda_4 +
\cdots } \\
a'(\lambda) &= a(\lambda') = (-1)^{\lambda'_2 + \lambda'_4+
\cdots}.\end{align*} We will set $a(\mu/\lambda) = a(\lambda)a(\mu)$
and similarly for $a'$. Note that $a(\mu/\lambda)$ does not depend
on $\mu$ and $\lambda$, but only depends on the set of squares which
lie in the difference of their Young diagrams.  If $\lambda \lessdot
\mu$ is a cover with a box added in the $i$-th row then
\begin{equation}
\label{eq:vsign} a'(\mu/\lambda) = (-1)^{\lambda_i} =
\begin{cases} 1 & \mbox{if $\mu/\lambda$ is a box on an odd
column,} \\ -1 & \mbox{if $\mu/\lambda$ is a box on an row column.}
\end{cases}
\end{equation}

Define the function $s_\alpha$ on covers $\lambda \lessdot \mu$ in
$Y$ by
$$
s_\alpha(\lambda \lessdot \mu) = (-1)^{\lambda_{1} + \lambda_{2} +
\cdots + \lambda_i}$$ if the box $\mu/\lambda$ is on the $i$-th row.
Define the function $s_\beta$ by $s_\beta(\lambda \lessdot \mu) =
a(\mu/\lambda)s_\alpha(\lambda \lessdot \mu)$. Thus we obtain two
signed posets $Y_\alpha = (Y,s_\alpha,a')$ and $Y_\beta =
(Y,s_\beta,a')$.

\psset{unit=0.7pt} \psset{linewidth=.5pt} \psset{framesep=2pt}
\begin{figure}[ht]
\pspicture(0,-10)(0,270)

\rput(0,0){\rnode{A}{$\emptyset$}}
\rput(0,50){\rnode{B}{\psframe(0,0)(10,10)}}

\rput(-70,100){\rnode{C1}{\psframe(0,0)(10,10)
\psframe(0,10)(10,20)}}
\rput(70,100){\rnode{C2}{\psframe[fillcolor=lightgray,fillstyle=solid](0,0)(10,10)
\psframe[fillcolor=lightgray,fillstyle=solid](10,0)(20,10)}}

\rput(-100,150){\rnode{D1}{\psframe(0,0)(10,10)
\psframe(0,10)(10,20)\psframe(0,20)(10,30)}}
\rput(0,150){\rnode{D2}{\psframe[fillcolor=lightgray,fillstyle=solid](0,0)(10,10)
\psframe[fillcolor=lightgray,fillstyle=solid](0,10)(10,20)
\psframe[fillcolor=lightgray,fillstyle=solid](10,10)(20,20)}}
\rput(100,150){\rnode{D3}{\psframe[fillcolor=lightgray,fillstyle=solid](0,0)(10,10)
\psframe[fillcolor=lightgray,fillstyle=solid](10,0)(20,10)
\psframe[fillcolor=lightgray,fillstyle=solid](20,0)(30,10)}}

\rput(-150,230){\rnode{E1}{\psframe(0,0)(10,10)
\psframe(0,10)(10,20)\psframe(0,20)(10,30) \psframe(0,30)(10,40)}}
\rput(-70,230){\rnode{E2}{\psframe[fillcolor=lightgray,fillstyle=solid](0,0)(10,10)
\psframe[fillcolor=lightgray,fillstyle=solid](0,10)(10,20)\psframe[fillcolor=lightgray,fillstyle=solid](0,20)(10,30)
\psframe[fillcolor=lightgray,fillstyle=solid](10,20)(20,30)}}
\rput(0,230){\rnode{E3}{\psframe(0,0)(10,10)
\psframe(10,0)(20,10)\psframe(0,10)(10,20) \psframe(10,10)(20,20)}}
\rput(70,230){\rnode{E4}{\psframe[fillcolor=lightgray,fillstyle=solid](0,0)(10,10)
\psframe[fillcolor=lightgray,fillstyle=solid](0,10)(10,20)
\psframe[fillcolor=lightgray,fillstyle=solid](10,10)(20,20)
\psframe[fillcolor=lightgray,fillstyle=solid](20,10)(30,20)}}
\rput(150,230){\rnode{E5}{\psframe(0,0)(10,10) \psframe(10,0)(20,10)
\psframe(20,0)(30,10) \psframe(30,0)(40,10)}}

\ncline{-}{A}{B} \mput*{$+$} \ncline{-}{B}{C1} \mput*{$-$}
\ncline{-}{B}{C2} \mput*{$-$} \ncline{-}{C1}{D1} \mput*{$+$}
\ncline{-}{C1}{D2} \mput*{$-$} \ncline{-}{C2}{D2} \mput*{$+$}
\ncline{-}{C2}{D3} \mput*{$+$}

\ncline{-}{D1}{E1} \mput*{$-$} \ncline{-}{D1}{E2} \mput*{$-$}
\ncline{-}{D2}{E2} \mput*{$-$} \ncline{-}{D2}{E3} \mput*{$-$}
\ncline{-}{D2}{E4} \mput*{$+$} \ncline{-}{D3}{E4} \mput*{$-$}
\ncline{-}{D3}{E5} \mput*{$-$}

\endpspicture
\caption{The signed poset $Y_\alpha$.  The shapes $\lambda$ such
that $a'(\lambda) = v(\lambda) = -1$ are shaded.} \label{fig:Yalpha}
\end{figure}

\begin{figure}[ht]
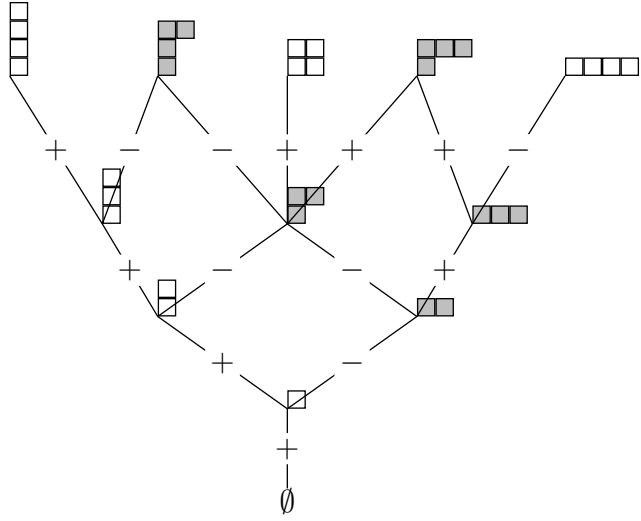

\pspicture(0,-10)(0,270)

\rput(0,0){\rnode{A}{$\emptyset$}}
\rput(0,50){\rnode{B}{\psframe(0,0)(10,10)}}

\rput(-70,100){\rnode{C1}{\psframe(0,0)(10,10)
\psframe(0,10)(10,20)}}
\rput(70,100){\rnode{C2}{\psframe[fillcolor=lightgray,fillstyle=solid](0,0)(10,10)
\psframe[fillcolor=lightgray,fillstyle=solid](10,0)(20,10)}}

\rput(-100,150){\rnode{D1}{\psframe(0,0)(10,10)
\psframe(0,10)(10,20)\psframe(0,20)(10,30)}}
\rput(0,150){\rnode{D2}{\psframe[fillcolor=lightgray,fillstyle=solid](0,0)(10,10)
\psframe[fillcolor=lightgray,fillstyle=solid](0,10)(10,20)
\psframe[fillcolor=lightgray,fillstyle=solid](10,10)(20,20)}}
\rput(100,150){\rnode{D3}{\psframe[fillcolor=lightgray,fillstyle=solid](0,0)(10,10)
\psframe[fillcolor=lightgray,fillstyle=solid](10,0)(20,10)
\psframe[fillcolor=lightgray,fillstyle=solid](20,0)(30,10)}}

\rput(-150,230){\rnode{E1}{\psframe(0,0)(10,10)
\psframe(0,10)(10,20)\psframe(0,20)(10,30) \psframe(0,30)(10,40)}}
\rput(-70,230){\rnode{E2}{\psframe[fillcolor=lightgray,fillstyle=solid](0,0)(10,10)
\psframe[fillcolor=lightgray,fillstyle=solid](0,10)(10,20)\psframe[fillcolor=lightgray,fillstyle=solid](0,20)(10,30)
\psframe[fillcolor=lightgray,fillstyle=solid](10,20)(20,30)}}
\rput(0,230){\rnode{E3}{\psframe(0,0)(10,10)
\psframe(10,0)(20,10)\psframe(0,10)(10,20) \psframe(10,10)(20,20)}}
\rput(70,230){\rnode{E4}{\psframe[fillcolor=lightgray,fillstyle=solid](0,0)(10,10)
\psframe[fillcolor=lightgray,fillstyle=solid](0,10)(10,20)
\psframe[fillcolor=lightgray,fillstyle=solid](10,10)(20,20)
\psframe[fillcolor=lightgray,fillstyle=solid](20,10)(30,20)}}
\rput(150,230){\rnode{E5}{\psframe(0,0)(10,10) \psframe(10,0)(20,10)
\psframe(20,0)(30,10) \psframe(30,0)(40,10)}}

\ncline{-}{A}{B} \mput*{$+$} \ncline{-}{B}{C1} \mput*{$+$}
\ncline{-}{B}{C2} \mput*{$-$} \ncline{-}{C1}{D1} \mput*{$+$}
\ncline{-}{C1}{D2} \mput*{$-$} \ncline{-}{C2}{D2} \mput*{$-$}
\ncline{-}{C2}{D3} \mput*{$+$}

\ncline{-}{D1}{E1} \mput*{$+$} \ncline{-}{D1}{E2} \mput*{$-$}
\ncline{-}{D2}{E2} \mput*{$-$} \ncline{-}{D2}{E3} \mput*{$+$}
\ncline{-}{D2}{E4} \mput*{$+$} \ncline{-}{D3}{E4} \mput*{$+$}
\ncline{-}{D3}{E5} \mput*{$-$}

\endpspicture
\caption{The signed poset $Y_\beta$.  The shapes $\lambda$ such that
$a'(\lambda) = v(\lambda) = -1$ are shaded.} \label{fig:Ybeta}
\end{figure}

The following theorem is similar to a calculation made in~\cite{Sta}, with a different definition of $U$ and $D$.

\begin{thm}
\label{thm:Ysigned} The signed posets $Y_\alpha$ and $Y_\beta$ are
weakly signed differential.
\end{thm}
\begin{proof}
We first prove the theorem for $Y_\alpha$.  Let $\lambda$ and $\mu$
be two distinct partitions satisfying $n = |\lambda| = |\mu|$. If
$\ip{(UD + DU)\lambda,\mu} \neq 0$ then it must be the case that
$\lambda \cap \mu = \nu$ where $|\nu| = n-1$. Let $\rho = \lambda
\cup \mu$. Suppose (without loss of generality) $\lambda/\nu$ lies
on the $i$-th row and $\mu/\nu$ lies on the $j$-th row, where $i <
j$. Then
\begin{align*}\ip{UD\,\lambda,\mu} &= (-1)^{\nu_{1} + \cdots +
\nu_{i}}(-1)^{\nu_{1}+ \cdots + \nu_{j}} a'(\lambda/\nu)
\\&= (-1)^{\lambda_{i+1} + \cdots + \lambda_j} a'(\lambda/\nu)
\end{align*}
and
\begin{align*}
\ip{DU\, \lambda,\mu} &= (-1)^{\lambda_{1} + \cdots +
\lambda_{j}}(-1)^{\mu_{1}+ \cdots + \mu_{i}} a'(\rho/\mu)
\\&= (-1)^{1+ \lambda_{i+1} + \cdots + \lambda_j} a'(\rho/\mu)
= \ip{-UD \, \lambda,\mu}
\end{align*}
using the fact that $\mu_i = \lambda_i + 1$ and the equality
$a'(\rho/\mu) = a'(\lambda/\nu)$.

Now we check that $\ip{(UD + DU)\lambda, \lambda} = 1$. We have
$$
UD \, \lambda = \left(\sum_{\mu \lessdot \lambda}
a'(\lambda/\mu)\right) \lambda
$$
and
$$
DU \, \lambda = \left(\sum_{\nu \gtrdot \lambda} a'(\nu/\lambda)
\right) \lambda.
$$
Using (\ref{eq:vsign}) we may pair up each inner corner of $\lambda$
with the outer corner of $\lambda$ in the next column to obtain the
required identity.  The coefficient of 1 arises from the outer
corner $\nu \gtrdot \lambda$ of $\lambda$ in the first column, which
has coefficient $a'(\nu/\lambda) = 1$.

Now for $Y_\beta$, the calculation of $\ip{(UD + DU)\lambda,
\lambda}$ is identical, while the calculation of
$\ip{UD\lambda,\mu}$ and $\ip{DU\lambda,\mu}$ is modified by a
factor of $a(\lambda)a(\mu)$ throughout.
\end{proof}

Obviously the edge labelings $s_\alpha, s_\beta$ can be modified in
other ways to still obtain a weakly differential poset.

\begin{thm}
\label{thm:Y} The signed poset $Y_\alpha$ is $\alpha$-signed
differential and the signed poset $Y_\beta$ is $\beta$-signed
differential.
\end{thm}
\begin{proof}
After Theorem~\ref{thm:Ysigned}, we need only check the equations
(\ref{eq:alpha}) and (\ref{eq:beta}).  For the $\beta$ case, the
coefficient of $\lambda$ in $(D-U)\P$ is
\begin{equation}
\sum_{\nu \gtrdot \lambda}s_\beta(\lambda \lessdot
\nu)a'(\nu/\lambda) - \sum_{\mu \lessdot \lambda}s_\beta(\mu
\lessdot \lambda) .
\end{equation}
If $\nu/\lambda$ is an outer corner in row $i$ of $\lambda$ and
$\lambda/\mu$ is the inner corner in row $i - 1$ then $$s_\beta(\mu
\lessdot \lambda) = (-1)^{\mu_1 + \cdots + \mu_{i-1}}
a(\lambda/\mu)= -a(\lambda/\mu)\,(-1)^{\lambda_1 + \cdots +
\lambda_{i-1}}$$ and using (\ref{eq:vsign}),
$$a'(\nu/\lambda)s_\beta(\lambda \lessdot \nu) =
(-1)^{\lambda_i}(-1)^{\lambda_1 + \cdots + \lambda_i}a(\nu/\lambda)
=(-1)^{\lambda_1 + \cdots + \lambda_{i-1}}a(\nu/\lambda).$$  Since
$a(\lambda/\mu) = - a(\nu/\lambda)$ the contributions of these two
corners cancel out. Finally for the unique outer corner
$\nu/\lambda$ in the first row we obtain a coefficient of
$(-1)^{\lambda_1}a'(\nu/\lambda)a(\nu/\lambda) = 1$.

For the $\alpha$ case, we note that without the additional factors
$a(\nu/\lambda)$ and $a(\lambda/\mu)$, the contributions of the
paired corners $\nu/\lambda$ and $\lambda/\mu$ will still cancel out
if we calculate $(U+D)\P$ instead.

\end{proof}

\begin{remark}
\label{rem:unique} Theorem~\ref{thm:Y} allows one to show that the
polynomials $g_k^\alpha(z), g_k^\beta(z)$ in Lemmas~\ref{lem:alpha}
and~\ref{lem:beta} and the coefficients $c_{ij}(w)$ in the proof of
Theorem~\ref{thm:w} are uniquely defined.  This follows from the
fact that in $Y_\alpha$ or $Y_\beta$, the element $U^n \hat 0$
contains the partition $(n)$ with non-zero coefficient and so is
non-zero.  Similarly, $D^n (n)$ is a non-zero multiple of $\hat 0$.
Thus one can ``extract'' the coefficients of $g^\alpha_k(z),
g^\beta_k(z)$ and $c_{ij}(w)$ one by one.
\end{remark}

\subsection{Fibonacci poset and signed reflection extensions}
\label{sec:refl}
Let $(P,s,v)$ be a signed poset such that $P = \cup_{0 \leq i \leq
n} P_n$.  Suppose $P$ is $\alpha$- or $\beta$-signed differential up
to level $n-1$.  In other words $UD + DU = I$ when restricted to
$\hat K ({\cup_{0 \leq i \leq n-1} P_i})$, and we have $(U+D)\P$ or
$(U-D)\P$ equal to $\sum_{0 \leq i \leq n-1} \P_i$ modulo $KP_n$.

We will now construct a signed poset $P^+ = \cup_{0 \leq i \leq n+1}
P^+_i$ with one more level than $P$. The signed poset
$(P^+,s^+,v^+)$ will satisfy $P^+_i = P_i$ for $0 \leq i \leq n$ and
also the equalities $s^+|_P = s$ and $v^+|_P = v$. First let
$P^+_{n+1}$ consist of elements $x^+$ for each $x \in P_{n-1}$ and
elements $y^*$ for each $y \in P_n$.  We will add the cover
relations $y^* \gtrdot y$ for each $y \in P_n$ and $x^+ \gtrdot y$
if $y \gtrdot x$, for each $x \in P_{n-1}, y \in P_n$. We then
define
$$ v^+(y^*) = v(y), \ \ \ v^+(x^+) = -v(x)$$
and
$$
s^+(y \lessdot y^*) = 1, \ \ \ s^+(y \lessdot x^+) = \begin{cases}
-v^+(x^+)v^+(y) s(x \lessdot y) & \mbox{in the $\alpha$ case,} \\
v^+(x^+)v^+(y) s(x \lessdot y) & \mbox{in the $\beta$ case.}
\end{cases}
$$
Let us use the notation $P^{+t}$ defined inductively $P^{+t} =
(P^{+(t-1)})^+$.  The following proposition is a signed analogue
of~\cite[Proposition 6.1]{Sta88}.

\begin{prop} \label{prop:ref}  Suppose $(P = \cup_{0 \leq i \leq n} P_n,s,v)$ is
$\alpha$- or $\beta$-signed differential up to level $n-1$.  Then
$(P^+,s^+,v^+)$ is $\alpha$- or $\beta$-signed differential up to
level $n$.  Thus ${\rm lim}_{t \to \infty} P^{+t}$ is a $\alpha$- or
$\beta$-signed differential poset.
\end{prop}
\begin{proof}
The proof is a straightforward case-by-case analysis: the signed
contribution of each cover cancels out with the contribution from
the reflected cover.
\end{proof}

\begin{remark}
Just as in the case of differential posets, the construction
described in Proposition~\ref{prop:ref} allows one to describe
infinitely many non-isomorphic signed differential posets.  They are
obtained by applying the signed reflection extension to the first
$n$-levels of the $\alpha$- or $\beta$-signed Young lattices.
\end{remark}

Let $Q = ({\hat 0},s,v)$ be the one element signed poset with
$v(\hat 0) = 1$.  Let $F_\alpha = (F_\alpha,s_\alpha,v_\alpha)$ and
$F_\beta = (F_\beta,s_\beta,v_\beta)$ denote the $\alpha$- and
$\beta$-signed differential posets obtained by the construction
${\rm lim}_{t \to \infty} Q^{+t}$.  Note that $v_\alpha = v_\beta$.
 We call these the ($\alpha$ or $\beta$) signed Fibonacci
differential posets.

We now give a non recursive description of $F_\alpha$ and $F_\beta$.
Define the {\it Fibonacci differential poset} $F = \cup_{r \geq 0}
F_r$ by letting $F_r$ be the set of words in the letters $\{1,2\}$
such that the sum of the letters is equal to $r$.  The covering
relations $x \lessdot y$ in $F$ are of two forms:
\begin{quote}
(a) $x$ is obtained from $y$ by changing a $2$ to a 1, provided that
the only letters to the left of this 2 are also 2's; or
\newline
(b) $x$ is obtained from $y$ by deleting the first 1 occurring in
$y$.
\end{quote}
The Fibonacci differential poset $F$ was defined in~\cite{Sta88} and independently in~\cite{Fom88} where it was called the Young-Fibonacci Lattice.

The poset $F$ is the underlying poset of the signed posets
$F_\alpha$ and $F_\beta$.  The reflection extension can be described
explicitly as follows.  The word $x^+$ is obtained from $x$ by
prepending a 2. The word $y^*$ is obtained from $y$ by prepending a
1.  If $a(x)$ denotes the number of 2's in $x$ then define $v(x) =
(-1)^a$ and
\begin{align}
\label{eq:fibalpha}
s'_\alpha(x \lessdot y)&= \begin{cases}(-1)^{i+1} & \mbox{if $x$ is $y$ with the (first) 1 in the $i$-th place deleted,} \\
(-1)^{i+1} & \mbox{if $x$ is $y$ with a 2 changed to a 1 in the
$i$-th place.}
\end{cases}\\
\label{eq:fibbeta}
s'_\beta(x \lessdot y) &= \begin{cases} (-1)^{i+1} & \mbox{if $x$ is $y$ with the (first) 1 in the $i$-th place deleted,} \\
(-1)^i & \mbox{if $x$ is $y$ with a 2 changed to a 1 in the $i$-th
place.}
\end{cases}
\end{align}

\begin{prop}
The signed posets $(F,s'_\alpha,v)$ and $(F,s'_\beta,v)$ are
identical to $(F_\alpha,s_\alpha,v_\alpha)$ and
$(F_\beta,s_\beta,v_\beta)$ respectively.
\end{prop}
\begin{proof}
The fact that the underlying posets are equal is straightforward to
verify (see also~\cite{Sta88}).  The equality $v = v_\alpha =
v_\beta$ follows immediately from induction.

To show that $s'_\alpha = s_\alpha$ and $s'_\beta = s_\beta$ we
again proceed by induction, using the following observations.
First, clearly the definitions agree on covers of the form $y
\lessdot y^*$.

If $x$ is obtained from $y$ by changing a $2$ to a $1$ in the $i$-th
position then $y$ is obtained from $x^+ = 2x$ by deleting a $1$ in
the $(i+1)$-position.  In this case we have $s_\alpha(x \lessdot y)
= -s_\alpha(y \lessdot x^+)$ and $s_\beta(x \lessdot y) = s_\beta(y
\lessdot x^+)$.

If $x$ is obtained from $y$ by deleting a $1$ in the $i$-position,
then $y$ is obtained from $x^+ = 2x$ by changing a 2 to a 1 in the
$i$-th position.  In this case we have $s_\alpha(x \lessdot y) =
s_\alpha(y \lessdot x^+)$ and $s_\beta(x \lessdot y) = -s_\beta(y
\lessdot x^+)$.
\end{proof}

The weighted sum of chains $e_\alpha(x)$ and $e_\beta(x)$ for
$F_\alpha$ and $F_\beta$ can be calculated explicitly.  We say that a word $x \in F$ is
{\it domino-tileable} if every non-initial, maximal, consecutive subsequence of 1's in $x$
has even length.  For example $x = 1221121111$ is domino-tileable but $y = 11212$ is not.

\begin{thm}
Suppose $x \in F$.  Then
\begin{align*}
e_\alpha(x) &= \begin{cases} 1 & \ \ \ \  \mbox{if $x$ is domino-tileable,} \\
                        0 & \ \ \ \ \mbox{otherwise.} \end{cases} \\
e_\beta(x) &= \begin{cases} v(x) & \mbox{if $x$ is domino-tileable,} \\
                        0 & \mbox{otherwise.} \end{cases}
\end{align*}
\end{thm}
\begin{proof}
We proceed by induction.  The result is clearly true for $x = \hat 0$, the empty word.  Now let $x \in F$
be an arbitrary word, and suppose $x = 2^j 1 w$ for some word $w$.  Then in $F$, the word $x$ covers the set of words
$C^-(x)= \{2^{i}12^{j-i-1}1w\} \cup \{2^jw\}$ where $0 \leq i \leq j-1$.  We use the recursive formula
$$
e(x) = \sum_{y \in C^-(x)} s(y \lessdot x)\, e(y).
$$
Suppose that the formula is known for all $y < x$.  If $j = 0$ the formula
is immediate from $s(y \lessdot 1y=x) = 1$.  For $j \geq 1$, the only possibly domino-tileable $y \in C^-(x)$ are $y_1 = 12^{j-1}1w$, $y_2 = 2^{j-1}11w$
and $y_3 = 2^jw$.    If $x$ is domino-tileable then only $y_1$ is (if $j = 1$, then $y_1 = y_2$), and $e(x) = s(y_1 \lessdot x) \, e(y_1)$, which
agrees with the theorem, using equations (\ref{eq:fibalpha}) and (\ref{eq:fibbeta}).

Otherwise, $x$ is not domino-tileable.  If $j = 1$, then $y_1 =
y_2$.  The word $y_1$ is domino-tileable if and only if $y_3$ is and
using equations (\ref{eq:fibalpha}) and (\ref{eq:fibbeta}) we see
that these contributions to $e(x)$ cancel out.  Thus we may assume
that $j > 1$.  Suppose first that $w$ is domino-tileable.  In this
case only $y_2$ and $y_3$ are domino-tileable and again their
contributions to $e(x)$ cancel out.  If $w$ is not domino tileable,
then none of $y_1,y_2,y_3$ are domino-tileable, so again $e(x) = 0$.

Finally we consider the case $x = 2^j$, where $j > 0$.  In this case
$x$ is domino-tileable but covers only a single domino-tileable word
$y = 12^{j-1}$.  Again the stated result follows inductively.
\end{proof}

\section{Sign-imbalance}
\label{sec:sign} We indicate here how our results can be applied to
sign-imbalance. If $P$ is a poset then a bijection $\omega: P \to
[n]= \{1,2,\ldots,n\}$ a is called a {\it labeling} of $P$.  A {\it
linear extension} of $P$ is an order-preserving map $f: P \to [n]$.
Given a labeled poset $(P,\omega)$ and a linear extension $f$ of
$P$, we obtain a permutation $\pi(f) = \omega(f^{-1}(1))
\omega(f^{-1}(2)) \cdots \omega(f^{-1}(n))\in S_n$.  We denote the
set of linear extensions of $(P,\omega)$ by $\L(P,\omega)$.  The
{\it sign-imbalance} of $(P,\omega)$ is the sum $I_{P,\omega}\sum_{f
\in \L(P,\omega)} {\rm sign}(\pi(f))$.  Up to sign, $I_{P,\omega}$
only depends on $P$.  Sign-imbalance was first studied by
Ruskey~\cite{Rus}.

\medskip

We first note a general basic property of the sign-imbalance of any
poset $P$ (see~\cite{Sta}).  Let $P$ be a finite poset with minimum
element $\hat 0$. We say that $P$ is {\it domino-tileable} if we can
find an increasing chain of order ideals (called a {\it domino
tiling})
$$
D = (I_0 \subset I_1 \subset \cdots \subset I_r = P)
$$
where the set theoretic difference $I_i - I_{i-1}$ for $1 \leq i
\leq r$ is a chain consisting of two elements and $|I_0| = 1$ or $0$
(depending on whether $|P|$ is odd or even).  Note that each domino
tiling $D$ of $P$ gives rise to a linear extension $f_D$ of $P$,
which is the unique linear extension satisfying $f_D(I_{i-1}) <
f_D(I_i-I_{-1})$.  The following Lemma follows from a sign-reversing
involution argument (see~\cite{Lam,Sta}).

\begin{lem}
\label{lem:dom}
Let $P$ be a finite poset with minimum element $\hat 0$ and $\omega: P \to \{1,2,\ldots,n\}$ any labeling of $P$.
If $P$ is domino-tileable, then
$$
I_{P,\omega} = \sum_D {\rm sign}(\pi(f_D))
$$
where the summation is over the domino-tilings of $P$.  If $P$ is not domino-tileable then $I_{P,\omega} = 0$.
\end{lem}

\subsection{Sign-imbalance of partition shapes}
First we consider the case of Young's lattice $Y$.  Let $\lambda$ be
a partition and $T$ a standard Young tableau (SYT) of shape
$\lambda$ (see~\cite{EC2}).  We will always draw our partition and
tableaux in English notation.  Picking the reverse of the standard
labeling of the poset $P_\lambda$ corresponding to the Young diagram
of $\lambda$, we can define the sign imbalance $I_\lambda$
explicitly as follows. The {\it reading word} $r(T)$ (or more
precisely the reverse row reading word) is the permutation obtained
from $T$ by reading the entries of $T$ from right to left in each
row, starting with the bottom row and going up. The {\it sign}
$s(T)$ is the sign of $r(T)$ as a permutation. Then the sign
imbalance is given by
$$
I_\lambda = \sum_T s(T)
$$
where the summation is over all standard Young tableaux $T$ with
shape $\lambda$.  We omit the labeling of the poset $P_\lambda$ in
our notation.

\begin{figure}[!ht]
$$
\tableau{{1}&{2}&{5}&{7}&{8} \\{3}&{6}&{9} \\{4}}
$$
\caption{A tableau $T$ with shape $531$, reading word $r(T) =
496387521$ and sign $s(T) = 1$.}
\end{figure}

\begin{remark}
Our reading order is the reverse of the reading order usually used
to define sign-imbalance for partitions~\cite{Lam,Rei,Sjo1,Sta}.
The resulting sign-imbalances differ by a factor of $(-1)^{{n
\choose 2}}$, where $n = |\lambda|$.
\end{remark}

\smallskip
We can connect the sign imbalance of Young diagrams with signed
differential posets as follows.  For $\lambda \in Y$ denote by
$e_\alpha(\lambda)$ and $e_\beta(\lambda)$ the signed sums of chains
in the signed posets $Y_\alpha$ and $Y_\beta$ respectively.

\begin{prop}
\label{prop:sign} Let $\lambda \in Y$.  Then $e_\alpha(\lambda) =
I_\lambda$ and $e_\beta(\lambda) = a(\lambda)I_\lambda$.
\end{prop}
\begin{proof}
A standard Young tableau $T$ of shape $\lambda$ is simply a maximal
chain $\emptyset =\lambda^{(0)} \subset \lambda^{(1)} \subset \cdots
\subset \lambda^{(l)} = \lambda$ in $Y$.  For a cover
$\lambda^{(i-1)} \lessdot \lambda^{(i)}$ on the $r$-th row, the sum
$\lambda^{(i-1)}_1 + \cdots + \lambda^{(i-1)}_r$ is equal to the number
of letters less than $i$ appearing after $i$ in $r(T)$.
\end{proof}

As corollaries we obtain the following theorem.

\begin{thm}
\label{thm:signim} Suppose $n \geq 2$.  Then
$$
\sum_{\lambda \vdash n} a'(\lambda)\,I_\lambda^2 = 0, \ \ \
\sum_{\lambda \vdash n} a'(\lambda)\,I_\lambda = 0, \ \ \
\sum_{\lambda \vdash n} a(\lambda) a'(\lambda)\,I_\lambda =
2^{\lfloor n/2 \rfloor}.
$$
\end{thm}
Theorem~\ref{thm:signim} was earlier conjectured in~\cite{Sta} and
proved independently in~\cite{Lam,Rei,Sjo1}.
\begin{proof}
Using Proposition~\ref{prop:sign} and Theorem~\ref{thm:Y}, the
result follows immediately from applying
Theorems~\ref{thm:signnfact} and~\ref{thm:signnsum} to $Y_\alpha$
and $Y_\beta$.
\end{proof}

One can define $I_{\lambda/\mu}$ for the Young diagrams of skew
shapes in an analogous manner to $I_\lambda$, by using a fixed
reading order. For our purposes, we suppose that we have picked a
reading order so that $I_{\lambda/\mu} = \ip{U^n \mu,\lambda}$ in
the $\alpha$-case and $I_{\lambda/\mu} = a(\lambda/\mu)\,\ip{U^n
\mu,\lambda}$ in the $\beta$-case.  Here $n = |\lambda/\mu|$.  This
differs from using the (reverse row) reading word by a factor of
$(-1)^a$, where $a$ is equal to the number of pairs $(s,t)$ of
squares $s \in \lambda/\mu$ and $t \in \mu$ such that $t$ is either
higher than or on the same row and to the left of $s$.

The following result,
due originally to Sj\"{o}strand~\cite{Sjo2}, follows from
Lemma~\ref{lem:skew}.
\begin{thm}
\label{thm:skew} Let $\lambda$ be a fixed partition and $n \in \N$.
Then
$$
\sum_{\mu/\lambda \vdash n} a'(\mu)I_{\mu/\lambda}^2 = \begin{cases}
\sum_{\lambda/\nu \vdash n} a'(\mu)I_{\lambda/\nu}^2 & \mbox{if $n$ is even,} \\
\sum_{\lambda/\nu \vdash n-1} a'(\nu)I_{\lambda/\lambda}^2 -
\sum_{\lambda/\nu \vdash n} a'(\nu)I_{\lambda/\nu}^2 & \mbox{if $n$
is odd.}\end{cases}
$$
\end{thm}
In Theorem~\ref{thm:skew}, the coefficients $a'(\mu)$ and $a'(\nu)$
can be replaced with $a(\mu)$ and $a(\nu)$ by making a similar
change in the definition of $Y_\alpha$.  In this form,
Theorem~\ref{thm:skew} is exactly Theorem 4.4 of~\cite{Sjo2}.  As a
consequence of Theorem~\ref{thm:F}, we have the following result.

\begin{thm}
Let
$$F_k(t) = \sum_{|\mu/\lambda| = k} a'(\mu/\lambda)\,I_{\mu/\lambda}^2 \,t^{|\lambda|}.$$
Then
$$
F_k(t) = \begin{cases} \prod_{i=1}^\infty \frac{1}{1-t^i} & \mbox{if
$k =
0$,} \\ \frac{1}{1+t}\prod_{i=0}^\infty \frac{1}{1-t^i} & \mbox{if $k = 1$,} \\
0 & \mbox{if $k \geq 2$.} \end{cases}
$$
\end{thm}
We have used the fact that $\sum_{\lambda \in Y} t^{|\lambda|} =
\prod_i \frac{1}{1-t^i}$. We write down one more result explicitly
from Proposition~\ref{prop:Gbeta}.
\begin{thm}
\label{thm:GY}
Let $l \in \N$.  Define
$$G_{2l}(t) = \sum_{|\lambda/\mu| = 2l}a(\lambda/\mu)\,I_{\lambda/\mu} \, t^{|\mu|}.$$
Then
\begin{align*}
G_{2l}(t) &= \left(\frac{2}{1+t^2}\right)^l \sum_{\lambda \in Y} a'(\lambda)t^{|\lambda|} \\ &= \left(\frac{2}{1+t^2}\right)^l \prod_{i=0}^\infty \left(\frac{1}{(1-t^{4i+1})(1+t^{4i+2})(1+t^{4i+3})(1-t^{4i+4})}\right).
\end{align*}
\end{thm}
\begin{proof}
We have $G_{2l}(t) = G_{2l}(Y_\beta,t)$ and so the first equation follows from Proposition~\ref{prop:Gbeta}.  The explicit expression for $\sum_{\lambda \in Y} a'(\lambda)t^{|\lambda|}$ follows directly from the definition
of $a'(\lambda)$.
\end{proof}
The constant term of the identity in Theorem~\ref{thm:GY} recovers
the $2^{\lfloor n/2 \rfloor}$ identity of Theorem~\ref{thm:signim}.
One can deduce many other results concerning sign-imbalance from
signed differential posets.  We leave this translation of our other
results, such as Proposition~\ref{prop:Galpha} to the reader.  An
interpretation of Theorem~\ref{thm:w} would require defining the
sign of an {\it oscillating tableau}.  However, such a definition
does not seem completely natural in the setting of sign-imbalance.


\subsection{Fibonacci sign-imbalance}
We now define the {\it Fibonacci distributive lattice} $\Fib$
(\cite{Sta88}).  The poset $\Fib$ has the same set elements as the
Fibonacci differential poset $F$, the set of words in the letters
$\{1,2\}$.  However the cover relations are defined differently.  If
$x = x_1 x_2 \cdots x_r$ and $y = y_1 y_2 \cdots y_l$ then $x
\lessdot y$ if $r \leq l$ and $x_i \leq y_i$ for $1 \leq i \leq r$.
The poset $\Fib$ is a distributive lattice.  It is equal to the
lattice of order ideals in an infinite dual (upside-down) tree.  For
$x \in \Fib$, we let $T_x$ denote the corresponding dual tree.  A
chain from $\hat 0$ to $x$ in $\Fib$ is a linear extension of $T_x$
which can be expressed simply as a tableau $T$ of the form shown in
Figure~\ref{fig:fib}.  The column lengths, read from left to right,
give the word $x$.
\begin{figure}[!h]
\begin{align*}
\tableau{{1}&{2}&{4}&{5}&{6}&{8} \\ {3}& & {7}& & {9}}
\end{align*}
\caption{A Fibonacci tableau $T$ with shape $212112$ and reading
word $r(T) = 312745698$.} \label{fig:fib}
\end{figure}
The top row is required to be increasing, and each column is also
required to be increasing. The reading word $r(T)$ of such a {\it
Fibonacci tableaux} is obtained by reading the columns from bottom
to top, starting with the leftmost column. This reading order
defines a labeling $\omega_x$ of $T_x$ and the corresponding
sign-imbalance is
$$
I_x = I_{T_x,\omega_x}  = \sum_{T} s(T)
$$
where $s(T)$ is the sign of $r(T)$ and $T$ varies over all Fibonacci tableaux with shape $x$.

We observe that a word $x$ is domino-tileable in the notation of Section~\ref{sec:refl} if and only
if the tree $T_x$ is domino-tileable as a poset.  Also recall that we define $v(x)$ to be $(-1)^{a(x)}$ where $a(x)$ is
the number of 2's in $x$.

\begin{prop}
\label{prop:fibagree} Let $x \in \Fib$.  Then
$$
I_x = \begin{cases} v(x) & \mbox{if $x$ is domino-tileable,} \\
0 &\mbox{otherwise.} \end{cases}
$$
Thus $I_x = e(x)$ when $x$ is considered an element of $F_\beta$.
\end{prop}
\begin{proof}
The proposition follows nearly immediately from Lemma~\ref{lem:dom}.
When $T_x$ is domino-tileable, it has a unique domino tiling with
corresponding linear extension of the form
$$
\tableau{{1}&{2}&{4}&{6}&{7}&{8} \\ &{3}&{5}&&&{9}}
$$
Here the boxes occupied by $\{2,3\}$, $\{4,5\}$, $\{6,7\}$ and
$\{8,9\}$ form the dominos.  The reading word of such a linear
extension has exactly as many inversions as columns of length 2, and
so $I_x = v(x)$.
\end{proof}
Proposition~\ref{prop:fibagree} is another example of enumerative
properties agreeing for the Fibonacci differential poset $F$ and the
Fibonacci distributive lattice $\Fib$ (see~\cite{Sta90}).  The
Fibonacci analogue of Theorem~\ref{thm:signim} is the following
result.
\begin{thm}
\label{thm:fibsignim} Suppose $n \geq 2$.  Then
$$
\sum_{x \in \Fib_n} v(x)\, I_x^2 = 0, \ \ \
\sum_{x \in \Fib_n } v(x)\, I_x = 2^{\lfloor n/2 \rfloor}.
$$
\end{thm}
\begin{proof}
Since $I_x = e(x)$ with $x$ considered an element of $F_\beta$, the result follows from applying
Theorems~\ref{thm:signnfact} and~\ref{thm:signnsum} to the $\beta$-signed differential poset $F_\beta$.
\end{proof}

Theorem~\ref{thm:fibsignim} is not difficult to prove directly.  For
example, it is easy to see that there are $2^{\lfloor n/2 \rfloor}$
domino-tileable $\{1,2\}$-words with sum $n$: this corresponds to
the $\lfloor n/2 \rfloor$ choices between a `2' or a `11', or
alternatively between a vertical domino and a horizontal domino. One
can obtain many more identities for Fibonacci sign-imbalance using
our results, and we leave the experimentation to the reader.  We
note the signed rank generating function
$$\sum_{x \in \Fib} v(x)\,t^{|x|} = \frac{1}{1-t+t^2}.$$


\begin{thebibliography}{a-a}
\bibitem{Fom88} {\sc S.~Fomin:} Generalized Robinson-Schensted-Knuth correspondence,
{\sl J. Soviet Math.} \textbf{41} (1988), 979-991.

\bibitem{Fom1} {\sc S.~Fomin:} Duality of graded graphs,
{\sl J. Algebraic Combin.} \textbf{3} (1994), 357--404.

\bibitem{Fom2} {\sc S.~Fomin:} Schensted algorithms for dual graded graphs,
{\sl J. Algebraic Combin.} \textbf{4} (1995), 5--45.

\bibitem{Lam} {\sc T.~Lam:} Growth diagrams, domino insertion and
sign-imbalance, {\sl J. Combin. Theory Ser. A} \textbf{107} (2004),
87--115.

\bibitem{Rei} {\sc A.~Reifergerste:} Permutation sign under the
Robinson-Schensted-Knuth correspondence, {\sl Ann. Combin.}
\textbf{8} (2004), 103--112.

\bibitem{Rus} {\sc F.~Ruskey:} Generating linear extensions of
posets by transpositions, {\sl J. Combin. Theory Ser. B} \textbf{54}
(1992), 77--101.

\bibitem{Sjo1} {\sc J.~Sj\"{o}strand:} On the sign-imbalance of
partition shapes, {\sl J. Combin. Theory Ser. A} \textbf{111}
(2005), 190--203.

\bibitem{Sjo2} {\sc J.~Sj\"{o}strand:} On the sign-imbalance of skew partition
shapes, preprint, 2005; {\tt math.CO/0507338}.

\bibitem{Sta88} {\sc R.~Stanley:} Differential posets, {\sl J. Amer.
Math. Soc.} \textbf{1} (1988), 919--961.

\bibitem{Sta90} {\sc R.~Stanley:}
Further combinatorial properties of two Fibonacci lattices, {\sl
Europ. J. Combin.} \textbf{11} (1990), 181--188.

\bibitem{Sta} {\sc R.~Stanley:} Some Remarks on Sign-Balanced
and Maj-Balanced Posets, {\sl Adv. Applied Math.} \textbf{34}
(2005), 88--902.

\bibitem{EC2} {\sc R.~Stanley}, {\sl Enumerative
Combinatorics, Vol 2}, Cambridge, 1999.

\bibitem{Var} {\sc A.~Varvak:}
Rook numbers and the normal ordering problem, {\sl J. Comb. Theory
Ser. A} \textbf{112} (2005), 292--307.

\bibitem{Whi} {\sc D.~White:}
Sign-balanced posets, {\sl J. Combin. Theory Ser. A.} \textbf{95}
(2001), 1--38.

\end{thebibliography}
\end{document}